\documentclass[12pt,oneside,english]{amsart}
\usepackage{lmodern}
\usepackage[T1]{fontenc}
\usepackage[latin9]{inputenc}
\usepackage{geometry}
\usepackage{color}
\usepackage{mathrsfs}
\usepackage{amsbsy}
\usepackage{amstext}
\usepackage{amsthm}
\usepackage{amssymb}
\usepackage[normalem]{ulem}
\usepackage{dsfont}
\usepackage[all,knot]{xy}
\xyoption{arc}
\usepackage{multirow}

\def\R{\mathbb R}
\def\C{\mathbb C}
\def\N{\mathbb N}
\def\Z{\mathbb Z}

\def\P{\mathbb P}

\def\fg{\mathfrak{g}}

\def\g{\boldsymbol{\mathfrak{g}}}

\def\s{\boldsymbol{\mathfrak{s}}}

\def\b{\boldsymbol{\mathfrak{b}}}

\def\BK{\boldsymbol{K}}
\def\BX{\boldsymbol{X}}
\def\l{\boldsymbol{\mathfrak{l}}}

\def\U{\mathcal{U}}

\makeatletter
\numberwithin{equation}{section}
\numberwithin{figure}{section}
\theoremstyle{plain}
\newtheorem{thm}{\protect\theoremname}[section]
\theoremstyle{remark}
\newtheorem{rem}[thm]{\protect\remarkname}
\theoremstyle{definition}
\newtheorem{defn}[thm]{\protect\definitionname}
\theoremstyle{plain}
\newtheorem{lem}[thm]{\protect\lemmaname}
\theoremstyle{plain}
\newtheorem{cor}[thm]{\protect\corollaryname}
\theoremstyle{plain}

\theoremstyle{plain}
\newtheorem{example}[thm]{\protect\examplename}

\makeatother

\usepackage{babel}
\providecommand{\corollaryname}{Corollary}
\providecommand{\definitionname}{Definition}
\providecommand{\lemmaname}{Lemma}
\providecommand{\remarkname}{Remark}
\providecommand{\theoremname}{Theorem}
\providecommand{\propositionname}{Proposition}
\providecommand{\examplename}{Example}

\begin{document}
\global\long\def\Crt{{\rm Crt}}%
\global\long\def\dist{{\rm \mathsf{dist}}}%
\global\long\def\Es{E_{\star}}%
\global\long\def\GS{\mathsf{GS}}%
\global\long\def\Rs{R_{\star}}%
\global\long\def\hRs{R_{\star}}%
\global\long\def\V{{\rm Vol}}%
\global\long\def\bs{\boldsymbol{\sigma}}%
\global\long\def\bsigma{\boldsymbol{\sigma}}%
\global\long\def\bx{\mathbf{x}}%
\global\long\def\BB{\mathbf{B}}%
\global\long\def\be{\mathbf{e}}%
\global\long\def\bxs{\mathbf{x}_{\star}}%
\global\long\def\qs{q_{\star}}%
\global\long\def\bz{\mathbf{z}}%
\global\long\def\by{\mathbf{y}}%
\global\long\def\sfv{{\mathsf{v}}}%
\global\long\def\sfu{{\mathsf{u}}}%
\global\long\def\bv{\mathbf{v}}%
\global\long\def\bu{\mathbf{u}}%
\global\long\def\bw{\mathbf{w}}%
\global\long\def\bn{\mathbf{n}}%
\global\long\def\bM{\mathbf{M}}%
\global\long\def\bJ{\mathbf{J}}%
\global\long\def\bG{\mathbf{G}}%
\global\long\def\indic{\mathds{1}}%
\global\long\def\diag{{\rm diag}}%
\global\long\def\sol{\mathsf{Sol}}%
\global\long\def\SN{\mathbb{S}^{N-1}}%
\global\long\def\BN{\mathbb{B}^{N}}%
\global\long\def\SNt{\mathbb{S}^{N-2}}%
\global\long\def\sp{{\rm sp}}%
\global\long\def\AA{\mathbb{A}}%
\global\long\def\E{\mathbb{E}}%
\global\long\def\P{\mathbb{P}}%
\global\long\def\R{\mathbb{R}}%
\global\long\def\V{\mathbb{V}}%
\global\long\def\T{\mathbb{T}}%
$ $
\global\long\def\BN{\mathbb{B}^{N}}%
\global\long\def\N{\mathsf{N}}%
\global\long\def\cpt{{\rm \mathsf{CP}}}%
\global\long\def\Cs{\mathscr{C}}%
\global\long\def\DD{\mathscr{D}}
\global\long\def\calE{\mathcal{E}}%
 
\global\long\def\epsilon{\varepsilon} 
\global\long\def\b{\beta}%
\global\long\def\BJ{\mathbf{J}}%
\global\long\def\ba{\mathbf{a}}%
\global\long\def\bb{\mathbf{b}}%
\global\long\def\bc{\mathbf{c}}%
\global\long\def\bcs{\mathbf{c}_{\star}}%
\global\long\def\grad{\nabla_{\mathrm{sp}}}%
\global\long\def\gradt{\bar \nabla}%
\global\long\def\Hess{\nabla_{\mathrm{sp}}^{2}}%
\global\long\def\Hesst{\bar\nabla^{2}}%
\global\long\def\HessT{\nabla_{\scriptscriptstyle T}^{2}}%
\global\long\def\cF{\mathcal{F}}%
\global\long\def\xim{\xi_{(m)}}
\global\long\def\xin{\xi_{(n)}}
\global\long\def\bxim{\bar{\xi}_{(m)}}
\global\long\def\bxik{\bar{\xi}_{(k)}}

\title{
Extensions of the constant family of    Harish-Chandra pairs 
of $SL_2(\mathbb{R})$}
\author{Eyal   Subag}
\begin{abstract}
We study and classify algebraic families of Harish-Chandra pairs over the complex affine line and over the complex projective line with   generic fiber that is  isomorphic to  the  Harish-Chandra pair of $SL_2(\R)$.
\end{abstract}

\maketitle
\tableofcontents 
\section{Introduction}
The study of algebraic families of  Harish-Chandra pairs and modules was initiated in \cite{MR4130851,MR4123111}. The main motivation there was a rigorous mathematical framework for   certain  limits of Lie groups and their representations that are common in physics and known as \textit{contractions} e.g., see \cite{MR55352,MR2942592}. 
Many of the algebraic families of  Harish-Chandra pairs in the literature are built out of a real reductive Lie group $G(\R)$. Such a group has a maximal compact subgroup $K(\R)$. The pair $(\fg,K)$ where $\fg$ is the complexification of the Lie algebra of $G(\R)$ and $K$ the complexification of $K(\R)$, is the Harish-Chandra pair of  $G(\R)$. 
The simplest construction of  an algebraic family of Harish-Chandra pairs over a smooth complex algebraic variety $\BX$, out of $G(\R)$ is a constant family. Explicitly,   the pair $(\fg,K)/\BX:=(\mathcal{O}_{\BX}\otimes_{\C}\fg, \BX\times_{\operatorname{spec}(\C)}K) $ is a constant algebraic family of   Harish-Chandra  pairs, where  $\mathcal{O}_{\BX}$ is the sheaf of regular functions on $\BX$.

Families of modules for  $(\fg,K)/\BX$ naturally  appear in representation theory of real reductive  groups see e.g.,  \cite{van2009analytic,MR4146144,MR4389792}.

In \cite{MR4130851}, two constructions associating  with a real reductive Lie group $G(\R)$,    a non-constant   algebraic family of Harish-Chandra pairs over $\mathbb{A}_{\C}^1$ were given. The first family $(\g_c,\BK)$ is called \textit{the contraction family} while the second family $(\g_d,\BK)$ is called \textit{the deformation family}. Both families share the same constant family of groups $\BK=\mathbb{A}_{\C}^1\times_{\operatorname{spec}(\C)}K$ and their restriction to $\C^{\times}\subset \mathbb{A}_{\C}^1$ is isomorphic to the constant family $(\fg,K)/\C^{\times}$ associated with  $G(\R)$. These  families have interesting applications. 

Families of Harish-Chandra modules for the deformation family  were used to give  a new algebraic construction  for the Mackey-Higson bijection \cite{MR3942563}. Full proof for   existence of  the Mackey-Higson bijection was given in \cite{MR4400734}. The bijection has a deep connection with the Baum-Connes conjecture \cite{MR3989145,MR1292018,MR1648112,MR2391803,MR1957086}. 

The contraction family and its modules, in the case of $G(\R)=SO(3,1)$,  were shown to play a central  role  in the fundamental quantum mechanical system of the hydrogen atom \cite{MR4460278,MR3827131}.

It is natural to ask whether other similar families of Harish-Chandra pairs can be associated with a real reductive Lie group.  The purpose of this paper is to  classify and study all  families of Harish-Chandra pairs over $\mathbb{A}_{\C}^1$ and over $\mathbb{P}_{\C}^1$  whose restriction to $\C^{\times}$ is isomorphic to the  constant family $(\mathfrak{sl}_2(\C),SO(2,\C))/\C^{\times}$
 associated with  $SL_2(\R)$. Such families of Harish-Chandra pairs are called extensions of the constant family $(\mathfrak{sl}_2(\C),SO(2,\C))/\C^{\times}$. 
 
We shall now review the structure of the  paper and some of  our main results.
In Section \ref{sec2} we recall the definitions of all relevant algebraic families. In Section \ref{sec3} we classify all extensions of $(\mathfrak{sl}_2(\C),SO(2,\C))/\C^{\times}$ over the complex affine line. We prove that the space of isomorphism classes of such families is in bijection with $\mathbb{N}_0$.  
As explained in Section \ref{sec2}, replacing an extension $(\g,\BK)$ by the space of global sections of $\g$ and the fiber $K$ of $\BK$ does not lose any information. 
In our case $K=SO(2,\C)$ and $\g$ becomes a Lie algebra over the polynomial ring $\C[x]$. Irreducible algebraic representations of $SO(2,\C)$ are parameterized by $\Z$. We let $\g_m$ denote the $K$ isotopic piece of global sections of $\g$ corresponding to the irreducible representation with parameter $m\in \Z$.

\begin{thm}[Sec. \ref{sec3}, Thm. \ref{th1'}]
Let  $(\g,\BK)$ be 
an extension of $(\mathfrak{sl}_2(\mathbb{C}),SO(2,\C))/\C^{\times}$.
Then   $\BK$ is isomorphic to $\mathbb{A}^1(\C)\times_{\operatorname{Spec}(\C)}SO(2,\C)$   and there exists a basis $\{{X},{Y},{H}\}$ of the space of global sections of $\g$ and $n\in \mathbb{N}_0$, such that  ${H}\in\g_0$, $X\in \g_2$, $Y\in \g_{-2}$ and    
\begin{eqnarray}\nonumber
&&[{H},\ {X}]=2 {X} , [{H},{Y} ]=-2 {Y},[{X},{Y} ]=x^n{H}.
\end{eqnarray}
Moreover, the section ${H}$ is unique and  the sections ${X}$ and ${Y}$ are unique up to simultaneous conjugation by  $SO(2,\C)$.  
\end{thm}
Denoting the extension corresponding to $n$ by $(\g(n),SO(2,\C))$. We show that $(\g(0),SO(2,\C))$ is the constant family, $(\g(1),SO(2,\C))$ is the contraction family and $(\g(2),SO(2,\C))$ is the deformation family. We prove that $(\g(1),SO(2,\C))$ is universal in the sense of the next theorem. 
\begin{thm}[Sec. \ref{sec3}, Thm. \ref{th2}]
Any extension  of $(\mathfrak{sl}_2(\mathbb{C}),SO(2,\C))/\C^{\times}$ 
is isomorphic to a pullback of
$(\g(1),SO(2,\C))$.
\end{thm}

In Section \ref{sec3.5} we show that any extension has a canonical realization inside the constant family. More precisely,
for $\fg:=\mathfrak{sl}_2(\mathbb{C})$, we let $\fg_m$, denote the isotypic spaces with respect to $SO(2,\C)$ corresponding to the irreducible representation parameterized by $m\in \Z$. A global section $F$ of the constant family $\g(0)$ is the same thing as an algebraic function from $\mathbb{A}_{\C}^1$ into $\fg$. Such a function has a well defined $n$-th derivative $F^{(n)}$ for any $n\in \mathbb{N}_0$.   
For   $n\in \mathbb{N}$, we let  $ \l(n)$ be the subspace of   the constant family $\g(0)$ consisting of all sections whose derivatives at $0$ up to degree $n$ belong to   $\fg_0\oplus \fg_2$, that is  \[ \l(n):=\{F\in\g(0)|F^{(m)}(0)\in \fg_0\oplus \fg_2, \forall m\in \{0,1,...,n-1\} \}.\]

\begin{thm}[Sec. \ref{sec3.5}, Thm. \ref{th6}]
  For   $n\in \mathbb{N}$, 
 the pair $(\l(n),SO(2,\C) )$   is a sub pair of $(\g(0),SO(2,\C))$ that is isomorphic to $(\g(n),SO(2,\C))$.
\end{thm}
 We also study the space of morphisms between extensions and the structure of the center of the universal enveloping algebra of each $\g(n)$. In Section \ref{sec4}  we classify extensions over $\mathbb{P}_{\C}^1$.

The author
is grateful to Joseph Bernstein for useful remarks. The author thanks Spyridon Afentoulidis-Almpanis for his comments on an earlier version of the manuscript. 

\section{Algebraic families of Harish-Chandra pairs}\label{sec2}
In this section  we  recall the  definition of an algebraic family of Harish-Chandra pairs over  a complex algebraic variety as given in \cite{MR4130851}.

Let $\BX$ be a smooth complex algebraic variety. We denote its   structure sheaf of regular functions by $\mathcal{O}_{\BX}$. An algebraic family of complex Lie algebras $\g$  over $\BX$ is a locally free sheaf of $\mathcal{O}_{\BX}$-modules which is also a sheaf of Lie algebras over $\BX$. In particular $\g$ is equipped with  Lie brackets  
\begin{eqnarray}\nonumber
 &&   [\_,\_]:\g\times \g\longrightarrow \g, 
\end{eqnarray}
which is $\mathcal{O}_{\BX}$ bi-linear, skew-symmetric and satisfying the Jacobi identity.  An algebraic family of complex algebraic groups $\boldsymbol{G}$ over $\BX$
is a smooth morphism between smooth complex algebraic varieties 
\begin{eqnarray}\nonumber
 &&     \boldsymbol{G} \longrightarrow \BX, 
\end{eqnarray}
that carries the structure of a group scheme over $\BX$. For every  family of complex algebraic groups there is an associated $\boldsymbol{G}$-equivariant family $\operatorname{Lie}(\boldsymbol{G})$ of complex Lie algebras over $\BX$. 
An algebraic family of Harish-Chandra pairs over $\BX$ is a pair $(\g,\BK)$  with $\BK$ being an algebraic family of complex algebraic groups over $\BX$, $\g$ being a $\BK$-equivariant algebraic family of Lie algebras over $\BX$. Along with a $\BK$-equivariant morphism of Lie algebras 
\begin{eqnarray}\nonumber
 &&   j:  \operatorname{Lie}(\BK) \longrightarrow \g, 
\end{eqnarray}
such that the action of $\operatorname{Lie}(\BK)$ on $\g$ that is obtained as the differential of the action of $\BK$ on $\g$ coincides with the action obtained by  $j$ composed with the adjoint action of $\boldsymbol{G}$ on $\g$. 
\begin{example}
Let $(\fg,K)$ be a complex Harish-Chandra pair. Note that this is the same thing as an algebraic family of Harish-Chandra pairs over the trivial one point complex algebraic variety. Then for any smooth complex affine algebraic variety $(\BX,\mathcal{O}_{\BX})$, the pair $(\mathcal{O}_{\BX}\otimes_{\C}\fg, \BX\times_{\operatorname{Spec}(\C)}K )$ is an algebraic family of Harish-Chandra pairs with all of its fibers canonically isomorphic to $(\fg,K)$. The equivariant Lie algebra morphism $j:\operatorname{Lie}(\BX\times_{\operatorname{Spec}(\C)}K)=\mathcal{O}_{\BX}\otimes_{\C} \operatorname{Lie}(K)  \longrightarrow \g$ is given by $j(f\otimes Y)=f j_0(Y)$ for every $f\in \mathcal{O}_{\BX}$ and $Y\in \operatorname{Lie}(K)$, and where $j_0:\operatorname{Lie}(K)\longrightarrow \fg$ is the morphism that is given with $(\fg,K)$. 
\end{example}

A family of Harish-Chandra pairs as in the last example  is called a \textit{constant} family and we shall denote it by $(\fg,K)/\BX$. Similarly a family of groups of the form   $\BX\times_{\operatorname{Spec}(\C)}K$ is called   constant, and a   family of Lie algebras of the form $\mathcal{O}_{\BX}\otimes_{\C}\fg$ is called constant.   Note that a section of $\mathcal{O}_{\BX}\otimes_{\C}\fg$ over an open subset  $U\subset \BX$ is canonically identified with an algebraic function from $U$ into $\fg$.

\subsection{Families of Harish-Chandra pairs over affine base}\label{subs1.1}
In this section we describe families of Harish-Chandra pairs with a constant group scheme and over affine base
 in terms of modules over a commutative ring.  This setup appears in many places throughout the paper.  
 
Assume that $\BX$ is a smooth complex affine algebraic variety. Denote its ring of global regular functions by $R_{\BX}$. Let  $(\g,\BK)$ be an algebraic family of Harish-Chandra pairs over $\BX$ with $\BK=\BX\times_{\operatorname{Spec}(\C)}K$ for some complex algebraic group $K$.  By looking on its space of global sections, such a family is equivalent to following data:
\begin{enumerate}
    \item A locally free  $R_{\BX}$-module $\g$ that is also a Lie algebra over $R_{\BX}$.(Note that by abuse of notation we denote the family of Lie algebras and its space of global sections by the same symbol). 
    \item An algebraic action 
    $\alpha:K\longrightarrow \operatorname{Aut}_{R_{\BX}}(\g)$    
    of $K$ on $\g$ by automorphisms of Lie algebras over $R_{\BX}$.
    \item A $K$-equivariant embedding of complex Lie algebras $j:\operatorname{Lie}(K)\longrightarrow \g$.
\end{enumerate}
Such that for every $Y\in \g$, $Z\in  \operatorname{Lie}(K)$,
\begin{eqnarray}\nonumber
 &&    [j(Z),Y]=d\alpha(Z)Y,
\end{eqnarray}
where $d\alpha$ stands for the action  of $\operatorname{Lie}(K)$ on $\g$ via endomorphisms 
 of Lie algebras over $R_{\BX}$, $d\alpha:\operatorname{Lie}(K)\longrightarrow \operatorname{End}_{R_{\BX}}(\g)$,  that is the differential of $\alpha$.
 
 Hence under the above mentioned assumptions
 we replace the pair  $(\g,\BK)$   
 by $(\g,K)$ where in the latter pair $\g$ is a Lie algebra over $R_{\BX}$ and $K$ a complex algebraic group.

Keeping the above assumptions of  $\BX$ being smooth affine and $\BK$ being constant,  if $(\g_1,\BX\times_{\operatorname{Spec}(\C)}K_1)$ and $(\g_2,\BX\times_{\operatorname{Spec}(\C)}K_2)$ are two algebraic families of Harish-Chandra pairs over $\BX$ with $K_1,K_2$ complex algebraic groups, then a morphism of pairs from  $(\g_1,\BX\times_{\operatorname{Spec}(\C)}K_1)$ to $(\g_2,\BX\times_{\operatorname{Spec}(\C)}K_2)$ is given by a pair of maps $(\psi_{\fg},\psi_K)$ with $\psi_{\fg}$ being a morphism of Lie algebras over $R_{\BX}$ and $\psi_K:K_1\longrightarrow K_2$ being a morphism of complex algebraic groups such that the following two conditions hold. 
\begin{enumerate}
    \item  \textit{Twisted equivariance}: $\psi_{\fg}(\alpha_1(k)X)=\alpha_2 (\psi_K(k)) \psi_{\fg}(X),$ $\forall X\in \g_1, k\in K_1$, and where $\alpha_i:K_i\longrightarrow \operatorname{Aut}_{R_{\BX}}(\g_i)$ is the action map for $i\in \{1,2\}$. 
    \item \textit{Compatibility with embeddings}: $\psi_{\fg}(j_1(Z))=j_2(d\psi_K(Z))$, $\forall Z\in \operatorname{Lie}(K_1)$ and where $j_i:\operatorname{Lie}(K_i)\longrightarrow \g_i$ is the $K_i$-equivariant embedding for $i\in\{1,2\}$, and $d\psi_K:\operatorname{Lie}(K_1)\longrightarrow \operatorname{Lie}(K_2)$ the differential of $\psi_K$.
\end{enumerate}
We denote the space of such morphisms by 
\[\operatorname{Hom}_{\text{pairs}}((\g_1,K_1),(\g_2,K_2)).\]
\begin{rem}\label{nr}
    When $\psi_{\fg}$ is an isomorphism and the center of $\g_1$ is trivial then compatibility with embeddings follows from twisted equivariance. 
\end{rem}
\begin{rem}\label{RemHom}
    For $(\psi_{\fg},\psi_K)\in \operatorname{Hom}_{\text{pairs}}((\g_1,K_1),(\g_2,K_2))$ with 
a one-to-one differential $d\psi_K:\operatorname{Lie}(K_1)\longrightarrow \operatorname{Lie}(K_2)$,
 we define a pullback pair $\psi_K^*(\g_2,K_2)$ of   $(\g_2,K_2)$ with respect to $\psi_K$, in which $\psi_K^*(\g_2,K_2)=(\g_2,K_1)$ and the action of $K_1$ on $\g_2$ is given by $\alpha_2\circ \psi_K:K_1\longrightarrow \operatorname{Aut}_{R_{\BX}}(\g_2)$, and the $K_1$-equivariant Lie algebras embedding is given by  
$j_2\circ d\psi_K :\operatorname{Lie}(K_1)\longrightarrow \g_2$.
Of course, if $(\psi_{\fg},\psi_K)\in \operatorname{Hom}_{\text{pairs}}((\g_1,K_1),(\g_2,K_2))$ then  $(\psi_{\fg},\mathbb{I})\in \operatorname{Hom}_{\text{pairs}}((\g_1,K_1),\psi_K^*(\g_2,K_2))$. 
\end{rem}
\section{Extensions of $(\mathfrak{sl}_2(\mathbb{C}),SO(2,\C))/\C^{\times}$ over   $\mathbb{A}^1_{\C}$}\label{sec3}
In this section we classify all extensions of $(\mathfrak{sl}_2(\mathbb{C}),SO(2,\C))/\C^{\times}$ over   $\mathbb{A}^1_{\C}$ and show that any extension is a pullback of the contraction family.

\subsection{Classification of  extensions  }\label{sec3.1}
Here we classify all extensions of $(\mathfrak{sl}_2(\mathbb{C}),SO(2,\C))/\C^{\times}$ over   $\mathbb{A}^1_{\C}$.  

Throughout the rest of the paper, we let $\BX$ be the complex algebraic variety $\mathbb{A}^1_{\C}$,  and let $\BX_0$ the subvariety consisting of $\BX$ take away a single point. Fixing a coordinate function on $\BX$ we shall identify $\BX$ with $\C$ and its ring of global regular functions with $\C[x]$, and $\BX_0$ with $\C^{\times}$  and its ring of global regular functions with $\C[x,x^{-1}]$.
 
\begin{defn}
 An  algebraic family of Harish-Chandra pairs $(\g,\BK)$ over $\mathbb{A}^1_{\C}$ is called an extension of $(\mathfrak{sl}_2(\mathbb{C}),SO(2,\C))/\C^{\times}$, if $\BK$ is a constant algebraic family of complex algebraic groups and 
 the restriction of $(\g,\BK)$ to $\C^{\times}$ is isomorphic to $(\mathfrak{sl}_2(\mathbb{C}),SO(2,\C))/\C^{\times}$.
\end{defn}

\begin{thm}\label{th1}
Let  $(\g,\BK)$ be an extension of $(\mathfrak{sl}_2(\mathbb{C}),SO(2,\C))/\C^{\times}$.
Then $\BK$ is isomorphic to $\mathbb{A}^1(\C)\times_{\operatorname{Spec}(\C)}SO(2,\C)$   and there exists a basis $\{{X},{Y},{H}\}$ of the space of global sections of $\g$ and $n\in \mathbb{N}_0$, such that  ${H}\in j(\operatorname{Lie}(\BK))$ and    
\begin{eqnarray}\nonumber
&&[{H},{X}]=2 {X} , [{H},{Y} ]=-2 {Y},[{X},{Y} ]=x^n{H}.
\end{eqnarray}
Moreover, the section ${H}$ is unique up to a sign and for a fixed ${H}$ the sections ${X}$ and ${Y}$ are unique up to simultaneous conjugation by  $SO(2,\C)$.  
\end{thm}
\begin{proof}
Since $\BK$ is a constant family  there exists a complex algebraic group $K$ such that $\BK=\BX\times_{\operatorname{Spec}(\C)}K$.
Let $\iota:\BX_0\longrightarrow \BX$ be the inclusion morphism. Since $\iota^*\BK\simeq \BX_0\times_{\operatorname{Spec}(\C)}K\simeq \BX_0\times_{\operatorname{Spec}(\C)}SO(2,\C)$, without loss of generality we can assume that $K=SO(2,\C)$.
Since $\BX$ and $\BX_0$ are smooth affine algebraic varieties and since all group schemes are constant we are in the setup of Subsection \ref{subs1.1} and we can consider the corresponding spaces of global sections. In this description of global sections over affine base, the   restriction $\iota^*\g=\g|_{\BX_0}$, of $\g$ to $\BX_0$ is given by the localization  $\C[x,x^{-1}]\otimes_{\C[x]}\g$. By our assumptions the family $(\g|_{\BX_0},SO(2,\C))$ is isomorphic to 
$(\mathfrak{sl}_2(\mathbb{C}),SO(2,\C))/\C^{\times}$.
We  realize  $(\mathfrak{sl}_2(\mathbb{C}),SO(2,\C))$ concretely as 
\begin{eqnarray}\nonumber
&& \mathfrak{sl}_2(\mathbb{C})=\left\{\left(\begin{matrix}
a & b\\
c & -a
\end{matrix}\right)\middle |a,b,c\in \mathbb{C}\right\}\\ \nonumber
&& SO(2,\C)=\left\{k(z):=\left(\begin{matrix}
z & 0\\
0 & z^{-1}
\end{matrix}\right)\middle |z\in \mathbb{C}^{\times}\right\}.
\end{eqnarray} 
We shall realize $\operatorname{Lie}(SO(2,\C))$ as
\begin{eqnarray}\nonumber
&& \mathfrak{so}_2(\mathbb{C})=\left\{\left(\begin{matrix}
a & 0\\
0 & -a
\end{matrix}\right)\middle |a\in \mathbb{C}\right\}.
\end{eqnarray} 
The embedding   $\tilde{j}_0:\operatorname{Lie}(SO(2,\C))\longrightarrow \C[x]\otimes_{\C}\mathfrak{sl}_2(\mathbb{C}) $
is given by $\tilde{j}_0(Z)=1\otimes j_0(Z)$ for every $Z\in \operatorname{Lie}(SO(2,\C))$, where $j_0:\operatorname{Lie}(SO(2,\C))\hookrightarrow \mathfrak{sl}_2(\mathbb{C})$ is the inclusion map. Any  isomorphism   \[(\psi_{\fg},\psi_K) \in   \operatorname{Hom}_{\text{pairs}}((\C[x,x^{-1}]\otimes_{\C}\mathfrak{sl}_2(\C),K),(\C[x,x^{-1}]\otimes_{\C[x]}\g,K))\] 
has exactly two possibilities for the 
algebraic isomorphisms  $\psi_K$ of $SO(2,\C)$. We denote them by $\psi_K^{\pm}$, 
explicitly 
\[\psi_K^{\pm }(g)=g^{\pm 1}, \quad \forall g\in K.\]
 We first assume that $\psi_K=\psi_K^{+}=\mathbb{I}_K$.
Since  $\g$ is a $K$ module (as a complex vector space) it has a canonical isotypical decomposition  $\g=\oplus_{n\in \mathbb{Z}}\g_n$ where for any $k(z)\in SO(2,\C)$ and $W\in \g_n$, we have $\alpha(k(z)) W=z^nW$. Since the actions of $\C[x]$ and  $K$   commute  each $\g_n$ is a $\C[x]$-module. Since $\C[x]$ is a principal ideal domain and since $\g$ is a free $\C[x]$-module then each  $\g_n$ is a free $\C[x]$-module. Similarly the isotypical decomposition of  $\C[x,x^{-1}]\otimes_{\C[x]}\g$ is given by $\C[x,x^{-1}]\otimes_{\C[x]}\g=\oplus_{n\in \mathbb{Z}}\left(\C[x,x^{-1}]\otimes_{\C[x]}\g_n\right)$. The twisted equivariance property and the fact that $\psi_K=\mathbb{I}_K$ imply that     $\psi_{\fg}$  is   an honest $K$-equivariant morphism  of Lie algebras over $\C[x,x^{-1}]$.
The $SO(2,\C)$ isotypic decomposition of  $\C[x,x^{-1}]\otimes_{\C}\mathfrak{sl}_2(\C)$  is given by 
\begin{eqnarray}\nonumber
&&\left(\C[x,x^{-1}]\otimes_{\C}e_{12}\right)\oplus  \left( \C[x,x^{-1}]\otimes_{\C}(e_{11}-e_{22})\right)\oplus 
    \left(\C[x,x^{-1}]\otimes_{\C}e_{21}\right) 
\end{eqnarray}
where $e_{ij}$ is the two-by-two matrix having a unique non-zero entry at the $(i,j)$-th entry which is equal to $1$, and  the action of $K$ on $\C[x,x^{-1}]\otimes_{\C}\mathfrak{sl}_2(\C)$ is determined  by  
\begin{eqnarray}\nonumber
&&k(z)\cdot 1\otimes (e_{11}-e_{22})=1\otimes (e_{11}-e_{22})\\ \nonumber
&& k(z)\cdot (1\otimes e_{12})=z^2\otimes e_{12}\\ \nonumber
&& k(z)\cdot (1\otimes e_{21})=z^{-2}\otimes e_{21}.
\end{eqnarray}
The chain of $K$-equivariant $\C[x]$-linear maps
\begin{eqnarray}\nonumber
    &&\g\longrightarrow \C[x,x^{-1}]\otimes_{\C[x]}\g  \longrightarrow \C[x,x^{-1}]\otimes_{\C}\mathfrak{sl}_2(\C) \\ \nonumber
    && s\longmapsto \quad \quad 1\otimes s \quad \quad \longmapsto \quad \quad (\psi_{\fg}^{\pm})^{-1}(1\otimes s)
\end{eqnarray}
takes a basis over $\C[x]$ for each isotypic piece to a basis over $\C[x,x^{-1}]$  of the corresponding isotypic piece. Hence $\g_n$ is a rank-one $\C[x]$-module for $n\in \{-2,0,2\}$ and zero otherwise.   For $n\in \{-2,0,2\}$, we let $X_n$ be a basis over $\C[x]$ of $\g_{n}$.
Since $\psi_{\fg}$ restricts to an isomorphism from any isotypic component of $\C[x,x^{-1}]\otimes_{\C}\mathfrak{sl}_2(\C)$ to the corresponding isotypic component of $\C[x,x^{-1}]\otimes_{\C[x]}\g $,  there must be $c_{-},c_0,c_{+}\in \C^{\times}$ and $m_{-},m_0,m_{+}\in \Z$ such that 
\begin{eqnarray}\nonumber
&& \psi_{\fg}(1\otimes e_{21})=c_{- }x^{m_{-}}\otimes X_{- 2} \\ \nonumber
&& \psi_{\fg}(1\otimes(e_{11}-e_{22}))=c_0x^{m_0}\otimes X_0 \\ \nonumber
&& \psi_{\fg}(1\otimes e_{12})=c_{+}x^{m_{+}}\otimes X_{ 2}.
\end{eqnarray} 

For  $i,j,k\in \{-2,0,2\}$, we let  $c^k_{i,j}(t)\in \mathbb{C}[x]$ be the structure constants of the Lie algebra $\g$ with respect to the ordered basis $\{X_{-2},X_0,X_2\}$. 
  Note that $c^k_{i,j}(x)$ are also the structure constants of $\C[x,x^{-1}]\otimes_{\C[x]}\g$ with respect to the ordered basis  $\{1\otimes X_{-2},1\otimes X_0,1\otimes X_2\}$. Clearly  
 \begin{eqnarray}\nonumber
 && [ \psi_{\fg}(1\otimes(e_{11}-e_{22})), \psi_{\fg}(1\otimes e_{12})]=\psi_{\fg}([1\otimes(e_{11}-e_{22}),1\otimes e_{12}]) 
 \end{eqnarray}
 which is equivalent to 
 \begin{eqnarray}\nonumber
 && [c_0x^{m_0}\otimes X_0,c_+x^{m_+}\otimes X_2 ]= \psi_{\fg}(2(1\otimes e_{12})=2c_+x^{m_+}\otimes X_2  \Longleftrightarrow  \\ \nonumber
 &&  \sum_{k\in \{ -2,0,2\}}c_0x^{m_0}c_+x^{m_+}c^k_{0,2}(x)\otimes {X}_k=2c_+x^{m_+}\otimes X_2.
  \end{eqnarray}  
 Hence   $c_{0,2}^{-2}(x)\equiv 0$, $c_{0,2}^0(x)\equiv 0$, and  $c_0x^{m_0}c_{0,2}^2=2$. As a result    $ c_{0,2}^2(x)= \frac{2}{c_0}x^{-m_0}$ and  $-m_0\in \mathbb{N}_0$.
Similarly, by examining the commutators $[ \psi_{\fg}(1\otimes(e_{11}-e_{22})), \psi_{\fg}(1\otimes e_{21})]$ and $[ \psi_{\fg}(1\otimes e_{12}), \psi_{\fg}(1\otimes e_{21})]$ we obtain that  $c_{0-,2}^0(x)\equiv0$, $c_{0,-2}^2(x)\equiv0$,  $ c_{0,-2}^{-2}(x)=- \frac{2}{c_0}x^{-m_0}$, $c_{2,-2}^{-2}(x)\equiv 0$, $c_{2,-2}^2(x)\equiv 0$ $c_{2,-2}^0(x)=\frac{c_0}{c_{-}c_+}x^{m_0-m_2-m_{-2}}$ and $m_0-m_2-m_{-2}\in \N_0$.

We define $n=m_0-m_+-m_{-}$. 
and ${H}=c_0X_0$, ${X}=c_{+}X_{ 2}$, ${Y}=c_{-}X_{-2}$,  and observe that
\begin{eqnarray}\nonumber
&&[{H},{X}]=2x^{-m_0}{X} , [{H},{Y} ]=-2x^{-m_0}{Y},[{X},{Y} ]=x^n{H}.
\end{eqnarray}
The embedding $1\otimes j$ of $\operatorname{Lie}(K)$ into $\C[x,x^{-1}]\otimes_{\C[x]}\g$ is given by 
\[(1\otimes j)(Z)=1\otimes j(Z), \quad \forall Z\in \operatorname{Lie}(K).\]
 By compatibility with embeddings of $(\psi_{\fg},\psi_K)$,   
\begin{eqnarray}\nonumber
     && (1\otimes j)\left(d\psi_K(e_{11}-e_{22}) \right)=\psi_{\fg}\left(\tilde{j}_0(e_{11}-e_{22})\right) \Longleftrightarrow \\ \nonumber
&& 1\otimes  j(e_{11}-e_{22})=\psi_{\fg}\left(1\otimes (e_{11}-e_{22})\right).
 \end{eqnarray}
 Hence $1\otimes  j(e_{11}-e_{22}) =c_0x^{m_0}\otimes X_0= x^{m_0}\otimes {H}$.
 This  equality is in $\C[x,x^{-1}]\otimes_{\C[x]}\g$. It implies that   $j(e_{11}-e_{22})= x^{m_0}{H}$ in $\g$. 
  Recall that $m_{0}\in -\mathbb{N}_0$ hence $m_0$ must be equal to $0$. Hence the commutation relations of the bases $\{ {H},{X}, {Y}\}$ are as stated in the theorem.  
   In addition ${H}= j(e_{11}-e_{22})\in 1\otimes \operatorname{Lie}(K)\subset  \operatorname{Lie}(\BK)$. 
   
  When   \[(\psi_{\fg},\psi_K) \in   \operatorname{Hom}_{\text{pairs}}((\C[x,x^{-1}]\otimes_{\C}\mathfrak{sl}_2(\C),K),(\C[x,x^{-1}]\otimes_{\C[x]}\g,K))\] 
  is an isomorphism in which 
  $\psi_K(g)=\psi_K^{-}(g)=g^{-1}$  for any $g\in K$, then as explained in Remark \ref{RemHom},
  \[(\psi_{\fg},\mathbb{I}_K) \in   \operatorname{Hom}_{\text{pairs}}((\C[x,x^{-1}]\otimes_{\C}\mathfrak{sl}_2(\C),K),\psi_K^*(\C[x,x^{-1}]\otimes_{\C[x]}\g,K)).\] 
  This allows us to use the previous case with 
   $\psi_K=\mathbb{I}_K$ and get the required basis.  In this case the basis $\{X,Y,H\}$ must satisfy $X\in \g_{-2}$, $Y\in \g_2$ and $H\in \g_0$ in contrary to case  with $\psi_K=\mathbb{I}_K$, in which  $X\in \g_{2}$, $Y\in \g_{-2}$ and $H\in \g_0$ . This finishes the proof of existence.

We shall now prove uniqueness.
 Given two bases with the required properties, $B:=\{  {H},{X}, {Y}\}$ and $B':=\{  {H}',{X}', {Y}'\}$, then the $\C[x]$-linear map $L:\g \longrightarrow \g$ satisfying $L({H})={H}'$, $L({X})={X}'$, $L({Y})={Y}'$ is an isomorphism of Lie algebras over  $\C[x]$.  
Both ${H}$ and ${H}'$ are   non zero complex multiple of $X_0$ and so there is $\alpha_H\in \C^{\times}$ such that $L({H})={H}'=\alpha_H {H}$. Since the eigenvalues of $\operatorname{ad}_{H}$ and of $\operatorname{ad}_{H'}$ must be equal to $\{-2,0,2\}$ we must have $\alpha_H=\pm 1$. If  $\alpha_H=1$, each of the elements ${X}$ and ${X}'$ is a basis for   $\g_2$ and  as a result the two are  proportional by an invertible element of $\C[x]$. Similarly each of   ${Y}$ and ${Y}'$ is a basis for $\g_{-2}$ and  as a result the two are  proportional by an invertible element of $\C[x]$. The commutation relations imply that the two proportionality constants are the inverse of each other. Namely, there exist $\lambda \in \C^{\times}=(\C[x])^{\times}$ such that  
\[{X}'=\lambda {X}, {Y}=\lambda^{-1} {Y}'. \]
For any square root of $\lambda$,   the matrix 
$k(\sqrt{\lambda}) \in SO(2,\C)$ satisfies
\[  k(\sqrt{\lambda}){X}=\lambda {X}={X}',  k(\sqrt{\lambda}){Y}=\lambda^{-1} {Y}={Y}'.  \] 
The case of $\alpha_H=-1$ is handled similarly. 
\end{proof}

\begin{rem}
By replacing $SO(2,\C)$ with $O(2,\C)$ in the above theorem, we can  show that any two bases with the required properties are in the same $O(2,\C)$-orbit.  
\end{rem}

\begin{cor}\label{c1}
    There is a bijection between $\mathbb{N}_0$ and  the set of equivalence classes of  
    extensions of $(\mathfrak{sl}_2(\mathbb{C}),SO(2,\C))/\C^{\times}$.
\end{cor}
\begin{proof}
We start by showing that if two extensions have bases as in the last theorem with the same $n\in \mathbb{N}_0$, then the extensions are isomorphic.
 As before we shall only consider the spaces of global sections for the families of Lie algebras and replace the constant group schemes with fiber $K=SO(2,\C)$ by their fiber. Let $(\g,K)$ and $(\g',K)$  be two extensions having bases  $B:=\{{X},{Y},{H}\}$ and $B':=\{{X'},{Y'},{H'}\}$ of the corresponding spaces of global sections  such that  ${H}\in j(\operatorname{Lie}(\BK))$, ${H'}\in j'(\operatorname{Lie}(\BK))$ and    
\begin{eqnarray} \nonumber
&&[{H}, {X}]=2 {X} , [{H},{Y} ]=-2 {Y},[{X},{Y} ]=x^n{H}, \\ \nonumber
&&[{H'}, {X'}]=2 {X'} , [{H'},{Y'} ]=-2 {Y'},[{X'},{Y'} ]=x^n{H'}.
\end{eqnarray}
 Define $\nu_{\fg}^{\pm}:\g\longrightarrow \g'$ via
\begin{eqnarray} \nonumber
&& \nu_{\fg}^{+}(H)= H', \hspace{2mm}\nu_{\fg}^{+}(X)= X', \hspace{2mm}\nu_{\fg}^{+}(Y)= Y',\\ \nonumber
    && \nu_{\fg}^{-}(H)= -H', \hspace{2mm}\nu_{\fg}^{-}(X)= Y', \hspace{2mm}\nu_{\fg}^{-}(Y)= X'
\end{eqnarray}
Both maps respect the Lie brackets, and exactly one of  $(\nu_{\fg}^{\pm},\mathbb{I}_K)$, is twisted-equivariance. By Remark, \ref{nr} it will also be compatible with the embeddings and hence an isomorphism of pairs. 

Now given $n\in \mathbb{N}_0$, we define  an algebraic family of Lie algebras  $\g$ via a basis $B:=\{{X},{Y},{H}\}$ for  its  space of global sections, namely $\g=\C[x]Y\oplus\C[x]H\oplus\C[x]X$. We define Lie brackets on $\g$ via   \begin{eqnarray}\nonumber
&&[{H}, {X}]=2 {X} , [{H},{Y} ]=-2 {Y},[{X},{Y} ]=x^n{H}.
\end{eqnarray} 
We define an action of $SO(2,\C)$ on $\g$ by declaring that the isotypic pieces are $\g_0=\C[x]H$, $\g_2=\C[x]X$, $\g_{-2}=\C[x]Y$.
The embedding $j:\mathfrak{so}(2,\C)\longrightarrow \g$ is the unique $\C$ linear map satisfying $j(e_{11}-e_{22})=H$. This is indeed an extension  of $(\mathfrak{sl}_2(\mathbb{C}),SO(2,\C))/\C^{\times}$.
\end{proof}

From the last proof we see that Theorem \ref{th1} can be reformulated as follows.

\begin{thm}\label{th1'}
Let  $(\g,\BK)$ be 
an extension of $(\mathfrak{sl}_2(\mathbb{C}),SO(2,\C))/\C^{\times}$.
Then   $\BK$ is isomorphic to $\mathbb{A}^1(\C)\times_{\operatorname{Spec}(\C)},SO(2,\C)$   and there exists a basis $\{{X},{Y},{H}\}$ of the space of global sections of $\g$ and $n\in \mathbb{N}_0$, such that  ${H}\in\g_0$, $X\in \g_2$, $Y\in \g_{-2}$ and    
\begin{eqnarray}\nonumber
&&[{H},\ {X}]=2 {X} , [{H},{Y} ]=-2 {Y},[{X},{Y} ]=x^n{H}.
\end{eqnarray}
Moreover, the section ${H}$ is unique and  the sections ${X}$ and ${Y}$ are unique up to simultaneous conjugation by an element of   $SO(2,\C)$.  
\end{thm}

For every $n\in \mathbb{N}_0$, we let $(\g(n),SO(2,\C))$ be a fixed representative for the algebraic family of Harish-Chandra pairs corresponding to $n$ from the last corollary. We fix a basis 
 $\{{X}_n,{Y}_n,{H}_n\}$ for $\g(n)$   such that  ${H}_n\in j(\operatorname{Lie}(SO(2,\C)))=\g(n)_0$, $X_n\in \g(n)_2$, $Y_n\in \g(n)_{-2}$ and    
\begin{eqnarray}\nonumber
&&[{H}_n,{X}_n]=2 {X}_n , [{H}_n,{Y}_n ]=-2 {Y}_n,[{X}_n,{Y}_n ]=x^n{H}_n.
\end{eqnarray}
For these choices we have $j(e_{11}-e_{22})=H_n$.
We shall keep these notations throughout the paper. 
\begin{rem}
    We note that the family $(\g(0),SO(2,\C))$ is the constant family $(\mathfrak{sl}_2(\mathbb{C}),SO(2,\C))/\C$, the family $(\g(1),SO(2,\C))$ is the \textit{contraction family} of $SL_2(\R)$, see section 2.1.3 of \cite{MR4130851}, and the family $(\g(2),SO(2,\C))$ is the \textit{deformation family} of $SL_2(\R)$, see section 2.1.2 of \cite{MR4130851}.
 \end{rem}

\subsection{The universality of $(\g(1),SO(2,\C))$}\label{21}
 
The next result  shows that $(\g(1),SO(2,\C))$ is universal among extensions of $(\mathfrak{sl}_2(\mathbb{C}),SO(2,\C))/\C^{\times}$.
\begin{thm}\label{th2}
Any extension  of $(\mathfrak{sl}_2(\mathbb{C}),SO(2,\C))/\C^{\times}$ 
is isomorphic to a pullback of
$(\g(1),SO(2,\C))$.
\end{thm}

\begin{proof}
Let $(\g,\BK)$ be an extension  of $(\mathfrak{sl}_2(\mathbb{C}),SO(2,\C))/\C^{\times}$. Recall that for any $\mu:\mathbb{A}_{\C}^1\longrightarrow \mathbb{A}_{\C}^1$, $\mu^*(\g)$ is the quotient of the $\C[x]$-module $\C[x]\otimes_{\C}\g $ by the submodule  generated by \[\{1\otimes f X-(f\circ \mu) \otimes X|f\in \C[x], X\in \g \}.\]
By abuse of notation the cosets in  $\mu^*(\g)$ shall be denoted by their representatives. The formula 
\[ [f\otimes X, g\otimes Y]=fg\otimes [X,Y], \quad f\otimes X, g\otimes Y\in \C[x]\otimes_{\C}\g \]
defines Lie brackets on $\mu^*(\g)$. 
 By corollary \ref{c1}, we can assume that there is some $n\in \mathbb{N}_0$, such that  $(\g,\BK)$ is equal to $(\g(n),SO(2,\C))$.
For $\mu_n(x)=x^n$, the $\C[x]$-linear map $L:\g(n)\longrightarrow \mu^*\g(1)$ satisfying
\begin{eqnarray}\nonumber
&&L(X_n)=1\otimes X_1,\hspace{2mm} L(Y_n)=1\otimes Y_1,\hspace{2mm} L(H_n)=1\otimes H_1
\end{eqnarray}
together with   $\psi_K^+:SO(2,\C)\longrightarrow SO(2,\C)$ defines  an isomorphism of families of Harish-Chandra pairs  from $(\g(n),SO(2,\C))$ onto $(\mu_n^*\g(1),SO(2,\C))$. 
\end{proof}
The next corollary shows that any extension is integrable, that is, it comes from an algebraic family of complex algebraic groups. 
\begin{cor}
Let $(\g,\BK)$ be 
an extension of $(\mathfrak{sl}_2(\mathbb{C}),SO(2,\C))/\C^{\times}$. Then 
 there is an algebraic family $\boldsymbol{G}$ of complex algebraic groups over  $\mathbb{A}_{\C}^1$ containing  $\BK$ and having its Lie algebra $\operatorname{Lie}(\boldsymbol{G})$ equal to $\g$.
\end{cor}
\begin{proof}
By the last Theorem we can assume that the extension $(\g,\BK)$ is equal to $(\mu_n^*\g(1),SO(2,\C))$  for some $n\in \mathbb{N}_0$. In \cite{MR3797197} it was shown that for the contraction family  $(\g_c,\BK)$ of a real reductive group $G(\R)$ there exists an algebraic family of complex algebra groups $\boldsymbol{G}_c$ over  $\mathbb{A}_{\C}^1$ containing $\BK$ and whose Lie algebra is equal to $\g_c$. We let $\boldsymbol{G}_1$ be the above mentioned family of groups in the case of $G(\R)=SL_2(\R)$. In this case $(\g_c,\BK)=(\g(1),SO(2,\C))$.  Since in our case the operation of pulling back with respect to $\mu_n$ and the operation of taking a Lie algebra for a family of groups commute, it follows that $\mu_n^*\boldsymbol{G}_1$  will have the required properties. 
\end{proof}

\begin{rem}
Using Theorem \ref{th2} we can give the space of equivalence classes of extensions the structure   of a monoid in which  multiplication is determined by composition of the maps $\mu_n$ defining the pullbacks.  With this monoid structure  the bijection in Corollary \ref{c1} is an isomorphism of monoids. 
\end{rem}

\section{Some properties of extensions  }\label{sec3.5}
In this sections we study the space of morphisms between any two extensions. We give a canonical embedding of any extension in the constant family and study the center of the universal enveloping algebra of any extension.
\subsection{Morphisms between extensions}  
In this section we study the space of morphisms between a given pair of extensions. 
For $m,n\in \mathbb{N}_0$, we  introduce the notations    
\begin{eqnarray}\nonumber
&&\operatorname{Hom}\nonumber(m,n):=\operatorname{Hom}_{\text{pairs}}((\g(m),SO(2,\C)),(\g(n),SO(2,\C))),\\ \nonumber
 &&\operatorname{Hom}(m,n)_{\pm}:=\{(\psi_{\fg},\psi_K)  
 \in \operatorname{Hom}(m,n)| \psi_K(k)=k^{\pm 1}, \forall k\in SO(2,\C)\}.
\end{eqnarray}
 Note that by remark \ref{RemHom}, we have a bijection betwenn $\operatorname{Hom}(m,n)_{+}$ and $\operatorname{Hom}(m,n)_{-}$.
We denote the space of complex polynomials in one variable of degree at most   $n$ by $\C_n[x]=\operatorname{span}_{\C}\{1,x,x^2,...,x^n\}$ and its subset of homogeneous polynomials by 
\[\C^{\text{hom}}_n[x]=\{f\in \C_n[x]|\exists k\in \mathbb{N}_0, f(tx)=t^kf(x), \forall t,x\in \C   \}.\]
\begin{thm}\label{t3}
The space $\operatorname{Hom}(m,n)$
is nonzero if and only if $m\geq n$. For $m\geq n$, 
\[  \operatorname{Hom}(m,n)_{+} \simeq  \C^{\text{hom}}_{m-n}[x] \simeq\operatorname{Hom}(m,n)_{-},\]
as sets.
\end{thm}

\begin{proof}
For $(\psi_{\fg},\psi_K)\in \operatorname{Hom}(m,n)_+$, we must have
\begin{eqnarray}\nonumber
&&\psi_{\fg}(X_m)=f_X(x) X_n,\hspace{2mm} \psi_{\fg}(Y_m)=f_Y(x) Y_n,\hspace{2mm} \psi_{\fg}(H_m)= f_H(x) H_n
\end{eqnarray}
for some $f_X,f_Y,f_H\in \C[x]$. Since $\psi_{\fg}$ is a morphism of Lie algebras, the commutation relations imply that 
\begin{eqnarray}\nonumber
&& f_H(x)f_X(x)=f_X(x), \hspace{1mm}  f_H(x)f_Y(x)=f_Y(x),\hspace{1mm}  f_X(x)f_Y(x)x^n=x^mf_H(x).
\end{eqnarray}
If $f_X(x)\equiv 0$, we must have $f_Y(x)\equiv 0$ and $f_H(x)\equiv 0$ and hence $\psi_{\fg}$ is the zero map. 
Similarly if $f_Y(x)\equiv 0$, $\psi_{\fg}$ is the zero map. Now assume that $f_X(x)f_Y(x)$ is not the zero polynomial then $f_H(x)\equiv 1$ and we must have $f_X(x)f_Y(x)x^n=x^m$. Hence we must have $n\leq  m$ and $f_X(x)f_Y(x)=x^{m-n}$. In particular the only possibilities  are 
\[f_X(x)=cx^k, f_Y(x)=c^{-1}x^{m-n-k}, \quad c\in \C^{\times}, k\in \{0,1,2,...,m-n \}. \]
We shall denote the above mentioned morphism by 
\[\psi_{\fg,m,n}^{c,k,+1}:\g(m)\longrightarrow \g(n) \]
The map sending $(\psi_{\fg},\psi_K)\in \operatorname{Hom}(m,n)_+$ to the unique $f_X(x)\in \C^{\text{hom}}_{m-n}[x]$ such that $\psi_{\fg}(X_m)=f_X(x)X_n$ is a bijection.
Similarly for $(\psi_{\fg},\psi_K)\in \operatorname{Hom}(m,n)_-$, for  $\psi_{\fg}$ to be nonzero we must have $n\leq m$ and  there exist $c\in \C^{\times}, k\in \{0,1,2,...,m-n \}$,  such that 
\begin{eqnarray}\nonumber
&&\psi_{\fg}(X_m)=cx^k Y_n,\hspace{2mm} \psi_{\fg}(Y_m)=c^{-1}x^{m-n-k} X_n,\hspace{2mm} \psi_{\fg}(H_m)=  -H_n.
\end{eqnarray}
We shall denote the above mentioned morphism by 
\[\psi_{\fg,m,n}^{c,k,-1}:\g(m)\longrightarrow \g(n) \]
The map sending $(\psi_{\fg},\psi_K)\in \operatorname{Hom}(m,n)_-$ to the unique $f_X(x)\in \C^{\text{hom}}_{m-n}[x]$ such that $\psi_{\fg}(X_m)=f_X(x)Y_n$ is a bijection.
\end{proof}

\begin{cor}
    Any extension  of $(\mathfrak{sl}_2(\mathbb{C}),SO(2,\C))/\C^{\times}$ can be embedded in the constant family $$(\mathfrak{sl}_2(\mathbb{C}),SO(2,\C))/\mathbb{A}^1_{\C}=(\g(0),SO(2,\C)).$$
\end{cor}
The following Lemma shall be needed later on. 
\begin{lem}\label{lem1}
  For any $m,n,p\in \mathbb{N}_0$, such that $m\geq n \geq p$,  any $k_1\in \{0,1,...,m-n\}$, $k_2\in \{0,1,...,n-p\}$, $c_1,c_2\in \C^{\times}$, $s_2,s_2\in \{-1,+1\}$,
    \[\psi_{\fg,m,n}^{c_1,k_1,s_1}\circ \psi_{\fg,n,p}^{c_2,k_2,s_2}=\psi_{\fg,m,p}^{c_1c_2,k_1+k_2,s_1s_2}\]  
\end{lem}
The proof is by direct calculation.

\subsection{Morphisms between localizations of extensions}  
In this section we study the space of morphisms between restrictions to $\C^{\times}$ of a given pair of extensions. This shall be needed in the last section of the paper where we classify extensions over $\mathbb{P}^1_{\C}$.

Recall  that $\iota^*\g(n)=\g(n)|_{\C^{\times}}=\C[x,x^{-1}]\otimes_{\C[x]}\g(n)$.  
For $m,n\in \mathbb{N}_0$, we define 
\[\operatorname{Hom}^0(m,n):=\operatorname{Hom}_{\text{pairs}}((\iota^*\g(m),SO(2,\C)),(\iota^*\g(n),SO(2,\C))).\]
\[\operatorname{Hom}^0(m,n)_{\pm}:=\{(\psi_{\fg},\psi_K)  
 \in \operatorname{Hom}^0(m,n)| \psi_K(k)=k^{\pm 1}, \forall k\in SO(2,\C)\}.\]
We let $\C((x))$ be the ring of formal Laurent series in one variable $x$ over $\C$, and let $\C((x))^{\text{hom}}$ be its subset of homogeneous elements. Explicitly, 
\[\C((x))^{\text{hom}}:=\{cx^n|n\in \Z, c\in \C \}.\]
\begin{thm}\label{th4}
For any $m,n\in \mathbb{N}_0$, there are bijections such that
\[  \operatorname{Hom}^0(m,n)_{+} \simeq  \C((x))^{\text{hom}} \simeq\operatorname{Hom}^0(m,n)_{-}.\]
\end{thm}
We omit the proof that  is similar to the proof of Theorem \ref{t3}.
For a later use we describe the hom spaces explicitly. For every $(\psi_{\fg},\psi_K)  
 \in \operatorname{Hom}^0(m,n)_{+}$ such that $\psi_{\fg}\neq 0$, we must have $\psi_{\fg}(1\otimes H_m)=  1\otimes H_n$, and
there exist $c\in \C^{\times}, k\in \mathbb{Z}$  such that
\begin{eqnarray}\nonumber
&&\psi_{\fg}(1\otimes X_m)=cx^k\otimes  X_n,\hspace{2mm} \psi_{\fg}(1\otimes Y_m)=c^{-1}x^{m-n-k} \otimes Y_n.
\end{eqnarray}
We shall denote the above mentioned morphism by 
\[\psi_{\fg^*,m,n}^{c,k,+1}:\iota^*\g(m)\longrightarrow \iota^*\g(n). \]
For every $(\psi_{\fg},\psi_K)  
 \in \operatorname{Hom}^0(m,n)_{-}$ such that $\psi_{\fg}\neq 0$, we must have $\psi_{\fg}(1\otimes H_m)=  -1\otimes H_n$, and
there exist $c\in \C^{\times}, k\in \mathbb{Z}$  such that
\begin{eqnarray}\nonumber
&&\psi_{\fg}(1\otimes X_m)=cx^k\otimes  Y_n,\hspace{2mm} \psi_{\fg}(1\otimes Y_m)=c^{-1}x^{m-n-k} \otimes X_n.
\end{eqnarray}
We shall denote the above mentioned morphism by 
\[\psi_{\fg^*,m,n}^{c,k,-1}:\iota^*\g(m)\longrightarrow \iota^*\g(n). \]

\subsection{Canonical realization of $(\g(n),SO(2,\C))$ inside the constant family}
In a previous section we saw that we can embed $(\g(n),SO(2,\C))$ in the constant family $(\g(0),SO(2,\C))$. However, this embedding is not unique and moreover  unless $n=0$ its image inside the constant family is not unique. 
In this section we show that for even $n$ there is a canonical realization of $(\g(n),SO(2,\C))$  inside the constant family. With an extra assumption  we  show that for any $n$ (even or odd) there is a  canonical realization of $(\g(n),SO(2,\C))$  inside the constant family. 

Recall that every $\sum_{i=1}^mf_i(x)\otimes Z_i\in \g(0)=\C[x]\otimes_{\C}\mathfrak{sl}_2(\mathbb{C})$ is canonically identified with the algebraic function from $\C$ into $\mathfrak{sl}_2(\mathbb{C})$ that is defined via
\[ t\longmapsto  \sum_{i=1}^mf_i(t)Z_i. \]
For such functions we have a well defined derivative of any order via differentiating the coefficients that are polynomial functions. We denote the $n$-th derivative of a function $F\in \g(0)$ by $F^{(n)}$. By convention $F^{(0)}=F$.
\subsubsection{The even case}
\begin{thm}
  For any $k\in \mathbb{N}$, the subset 
  \[ \s(2k):=\{F\in\g(0)|F^{(m)}(0)\in \mathfrak{so}(2,\C), \forall m\in \{0,1,...,k-1 \} \}\]
  is sub Lie algebra over $\C[x]$ that together with $SO(2,\C)$ is a sub pair of $(\g(0),SO(2,\C))$ that is isomorphic to $(\g(2k),SO(2,\C))$.
\end{thm}

\begin{proof}
We let $\{{X},{Y},{H}\}$ be
an ordered basis of $\mathfrak{sl}_2(\mathbb{C})$
   such that  ${H}\in  \mathfrak{so}(2,\C)$, $X\in (\mathfrak{sl}_2(\mathbb{C}))_2 $, $Y\in (\mathfrak{sl}_2(\mathbb{C}))_{-2} $ and    
\begin{eqnarray}\nonumber
&&[{H},{X}]=2 {X} , [{H},{Y} ]=-2 {Y},[{X},{Y} ]={H}.
\end{eqnarray}
The set 
$\{{X}_0:=1\otimes X,{Y}_0:=1\otimes Y,H_0:=1\otimes {H}\}$ is an ordered basis of $\g(0)$ with the same structure constants as $\{{X},{Y},{H}\}$. Note that the notation $\{{X}_0,{Y}_0,{H}_0\}$ is consistent with our fixed basis for $\g(0)$ form  Section \ref{sec3.1}.
 Clearly 
$B:=\{1\otimes H, x^k\otimes X,  x^k\otimes Y \}$ is linearly independent set over $\C[x]$ that belongs to  $\s(2k)$. Any element $F$ of  $\g(0)$ can be written uniquely in the form $f_X(x)\otimes X+f_Y(x)\otimes Y+f_H(x)\otimes H$, for some $f_X(x),f_Y(x),f_H(x)\in \C[x]$. For any $m\in \mathbb{N}_0$, 
\[F^{(m)}=f_X^{(m)}(x)\otimes X+f_Y^{(m)}(x)\otimes Y+f_H^{(m)}(x)\otimes H.\]
Hence $F\in \s(2k)$ if and only if
\[f_X^{(m)}(0)X+f_Y^{(m)}(0) Y+f_H^{(m)}(0) H\in \C H, \quad \forall m\in \{0,1,...,k-1\} .\]
This shows that $B$ is a basis for $\s(2k)$. Since  
\begin{eqnarray}\nonumber
    && [1\otimes H,x^k\otimes X]=2x^k\otimes X\\ \nonumber 
    && [1\otimes H,x^k\otimes Y]=-2x^k\otimes Y\\ \nonumber 
    && [x^k\otimes X,x^k\otimes Y]=x^{2k}\otimes H.
\end{eqnarray}
it follows that $(\s(2k),SO(2,\C))\simeq (\g(2k),SO(2,\C))$.
\end{proof}
\subsubsection{The general case}
For $n\in \mathbb{N}$, there are many different isomorphic copies of $\g(n)$ in $\g(0)$. For even $n$, as explained in the last theorem, there is a natural choice for such a copy, given by $\s(n)$. Here we shall make a  a different  choice that will work equally well for all $n\in \mathbb{N}$. This choice is natural for a given positive system of $\mathfrak{sl}_2(\mathbb{C})$.

As in the proof of the last Theorem, we let $\{{X},{Y},{H}\}$ be
an ordered basis of $\fg:=\mathfrak{sl}_2(\mathbb{C})$
   such that   ${H}\in  \mathfrak{so}(2,\C)=g_0$, $X\in g_2 $, $Y\in g_{-2} $ and    
\begin{eqnarray}\nonumber
&&[{H},{X}]=2 {X} , [{H},{Y} ]=-2 {Y},[{X},{Y} ]={H}.
\end{eqnarray}

\begin{thm}\label{th6}
  For any $n\in \mathbb{N}$, the subset 
  \[ \l(n):=\{F\in\g(0)|F^{(m)}(0)\in \fg_0\oplus \fg_2, \forall m\in \{0,1,...,n-1\} \}\]
  is sub Lie algebra over $\C[x]$ that together with $SO(2,\C)$ is a sub pair of $(\g(0),SO(2,\C))$ that is isomorphic to $(\g(n),SO(2,\C))$.
\end{thm}
\begin{proof}
    As in the proof of the last Theorem, we can show that the set 
$\{1\otimes X,x^n\otimes Y,1\otimes {H}\}$ is a basis over $\C[x]$ for $\l(n)$ satisfying 
\begin{eqnarray}\nonumber
    && [1\otimes H,1\otimes X]=2(1\otimes X)\\ \nonumber 
    && [1\otimes H,x^n\otimes Y]=-2(x^n\otimes Y)\\ \nonumber 
    && [1\otimes X,x^n\otimes Y]=x^{n}(1\otimes H).
\end{eqnarray}
This implies that $(\l(n),SO(2,\C)\simeq (\g(n),SO(2,\C))$.
\end{proof}

\subsection{The center  of $\mathcal{U}(\g(n))$}
In this section we explicitly describe the center of  the universal enveloping algebra of $\g(n)$. This has applications to representation theory of $\g(n)$ and of $(\g(n),SO(2,\C))$.

As above   we fix a basis $\{{X},{Y},{H}\}$  of $\fg:=\mathfrak{sl}_2(\mathbb{C})$
   such that  ${H}\in  \mathfrak{so}(2,\C)$ and    
\begin{eqnarray}\nonumber
&&[{H},{X}]=2 {X} , [{H},{Y} ]=-2 {Y},[{X},{Y} ]={H}.
\end{eqnarray}
It is well known that $\mathcal{Z}(\fg)$, the center of the universal enveloping algebra of $\fg$ is equal to the polynomial algebra $\C[\Omega_{\fg}]$, where $\Omega_{\fg}$ is the canonical Casimir element which is given in terms of the mentioned basis by 
\[\Omega_{\fg}= \frac{1}{8}H^2+\frac{1}{4}(XY+YX) .\]

In this section we shall realize the extensions $(\g(n),SO(2,\C))$ concretely inside the constant family $ (\g(0),SO(2,\C))$ via the isomorphism
$(\g(n),SO(2,\C))\simeq (\l(n),SO(2,\C)) \subset (\g(0),SO(2,\C)) $. We shall not explicitly mention the isomorphism and simply work with  $\l(n)$ and the basis 
\[B_n:=\{H_n:=1\otimes H, X_n:= 1\otimes X, Y_n:=x^n\otimes Y\}\subset \g(0)\] of $\l(n)$.
We define 
\[\Omega_{n}:= x^n\otimes \Omega_{\fg} \in \C[x]\otimes_{\C}\mathcal{Z}(\fg)\simeq \mathcal{Z}(\l(0)) .\]
\begin{thm}
   The center of the universal enveloping algebra of $\l(n)$ is equal to the commutative polynomial algebra over $\C[x]$ generated by $\Omega_{n}$ inside $\C[x]\otimes_{\C}\mathcal{Z}(\fg)\simeq \mathcal{Z}(\g(0))$.  In particular $\mathcal{Z}(\l(n))=\C[x,\Omega_{n}]$.
\end{thm}

\begin{proof}
    We first note that 
    \[\Omega_{n}= \frac{x^n}{8}H_n^2+\frac{1}{4}(X_nY_n+Y_nX_n)\in  \mathcal{Z}(\l(0)) .\]
    The proof is similar to the proof of Lemma 6.4 in  \cite{MR3942563}.
\end{proof}
 
\begin{rem}
The family of commutative $\C[x]$-algebras $\mathcal{Z}(\l(n))$ has a canonical subspace (which is a free rank one $\C[x]$-module) consisting of all elements of $\mathcal{Z}(\l(n))$ that are of degree two in the sense of PBW-filtration and act on the trivial $\mathcal{U}(\l(n))$-module $\C[x]$ by zero. This subspace, known as the Casimir subfamily of $\mathcal{Z}(\l(n))$, has $\Omega_{n}$ as its basis. The basis and hence $\Omega_{n}$ is unique up to a nonzero complex multiplicative constant.
\end{rem}

 \section{Extensions of $(\mathfrak{sl}_2(\mathbb{C}),SO(2,\C))/\C^{\times}$ over    $\mathbb{P}_{\C}^1$}\label{sec4}  

In this section we classify extensions of $(\mathfrak{sl}_2(\mathbb{C}),SO(2,\C))/\C^{\times}$ over    $\mathbb{P}_{\C}^1$.

\subsection{Extensions over  $\mathbb{P}_{\C}^1$}
Following tradition we realize $\mathbb{P}_{\C}^1$ via homogeneous coordinates as $\{[z_1;z_2]\in \C^2 \}/\sim $ where $[z_1;z_2] \sim [cz_1;cz_2]$ for any $c\neq 0$, $z_2,z_2\in \C$. We identify  $\C^{\times}$ with $\{[z_1;z_2]\in \mathbb{P}_{\C}^1| z_1z_2\neq 0 \}$.
\begin{defn}
 An  algebraic family of Harish-Chandra pairs $(\g,\BK)$ over $\mathbb{P}_{\C}^1$ is called an extension of $(\mathfrak{sl}_2(\mathbb{C}),SO(2,\C))/\C^{\times}$, if $\BK$ is a constant algebraic family of complex algebraic groups and 
 the restriction of $(\g,\BK)$ to $\C^{\times}$ is isomorphic to $(\mathfrak{sl}_2(\mathbb{C}),SO(2,\C))/\C^{\times}$.
\end{defn}

For $i\in \{1,2\}$, we define $U_i:=\{[z_1;z_2]\in\mathbb{P}_{\C}^1| z_i\neq 0 \}$ and  $\varphi_i:\U_i\longrightarrow \C$, via $\varphi_1([z_1;z_2])=z_2z_1^{-1}$,   
$\varphi_2([z_1;z_2])=z_1z_2^{-1}$. Then $U_1\cup U_2$ is an open affine cover of  $\mathbb{P}_{\C}^1$  and $\varphi_i$ are algebraic isomorphisms. Moreover $\varphi_2^{-1}\circ \varphi_1|_{U_1\cap U_2}:U_1\cap U_2\longrightarrow U_1\cap U_2$ is an isomorphism of affine algebraic varieties. It is  well known that a sheaf on $\mathbb{P}_{\C}^1$ is the same thing as a pair of sheaves, one on $U_1$ and the other on $U_2$, that agree on $U_1\cap U_2$.  Using the isomorphisms $\varphi_i$, a sheaf on $\mathbb{P}_{\C}^1$ is an ordered  pair of sheaves $(\mathcal{F}_1,\mathcal{F}_2)$ both on $\C$, together with an isomorphism
$\psi:\mathcal{F}_1|_{\C^{\times}}\longrightarrow \mathcal{F}_2|_{\C^{\times}}$ between their restrictions to $\C^{\times}$. We shall denote the corresponding sheaf on $\mathbb{P}_{\C}^1$ by $(\mathcal{F}_1,\mathcal{F}_2,\psi)$ and its isomorphism class by  $[(\mathcal{F}_1,\mathcal{F}_2,\psi)]$.

 \subsection{Classification of extensions over  $\mathbb{P}_{\C}^1$}

\begin{thm}
 The collection of all the equivalence classes 
\[  \left[\left((\g(m),K),  (\g(n),K), (\psi_{\fg^*,m,n}^{1,k,1},\psi_K^{1})\right)\right], \] 
as $m$ and $n$ vary in $\mathbb{N}_0$ and  $k$ vary in  $\Z$,  
is an exhaustive list of  equivalence classes of extensions  of $(\mathfrak{sl}_2(\mathbb{C}),SO(2,\C))/\C^{\times}$ over $\mathbb{P}_{\C}^1$. Moreover, the  representatives in the list are pairwise non-isomorphic.   
\end{thm}
\begin{proof}
From the above discussion together with Theorem \ref{th1} and  Theorem \ref{th4} any extension over $\mathbb{P}_{\C}^1$
is isomorphic to some \[  \left((\g(m),K),  (\g(n),K), (\psi_{\fg^*,m,n}^{c,k,\pm1},\psi_K^{\pm})\right).  \]  
Note  that the automorphism 
\[(\psi_{\fg,m,m}^{c^{-1},0,\pm1},\psi_K^{\pm1}):(\g(m),K)\longrightarrow (\g(m),K) \]
satisfies $ \psi_{\fg^*,m,m}^{c^{-1},0,\pm1}=1\otimes \psi_{\fg,m,m}^{c^{-1},0,\pm1}$ where 
\[(1\otimes \psi_{\fg,m,m}^{c^{-1},0,\pm1},\psi_K^{\pm1}):(\iota^* \g(m),K)\longrightarrow (\iota^*\g(m),K) \]
is the induced map  between  the corresponding localizations.
Since 
\[   (\psi_{\fg^*,m,n}^{c,k,\pm1},\psi_K^{\pm})\circ (\psi_{\fg^*,m,m}^{c^{-1},0,\pm1},\psi_K^{\pm1})= (\psi_{\fg^*,m,n}^{1,k,1},\psi_K^{1}),\]
we may assume that $c=1$ and $\psi_K$ is the identity map. This proves exhaustiveness. 
We now show that the  representatives in the above list are pairwise non-isomorphic.
Assume that for $i\in \{1,2\}$,
\[  \left((\g(m_i),K),  (\g(n_i),K), (\psi_{\fg^*,m_i,n_i}^{1,k_i,1},\psi_K^{1})\right)  \]
is an extension  of $(\mathfrak{sl}_2(\mathbb{C}),SO(2,\C))/\C^{\times}$ over $\mathbb{P}_{\C}^1$. 
Assume in addition that the two extensions are isomorphic. This means that $(\g(m_1),K)$ is isomorphic to $(\g(m_2),K)$ and $(\g(n_1),K)$ is isomorphic to $(\g(n_2),K)$ and the two isomorphisims are compatible with $(\psi_{\fg^*,m_1,n_1}^{1,k_1,1},\psi_K^{1})$ and $(\psi_{\fg^*,m_2,n_2}^{1,k_2,1},\psi_K^{1})$. Hence we must have $m_1=m_2$ and $n_1=n_2$. We denote the corresponding isomorphisms  by  $(\psi_{\fg,m_1,m_1}^{c_m,0,s_m},\psi_K^{s_m})$ and $(\psi_{\fg,n_1,n_1}^{c_n,0,s_n},\psi_K^{s_n})$. Compatibility means that the following diagram, is commutative

\[\xymatrix{ (\iota^*\g(m_1),K)\ar[dd]^{(1\otimes \psi_{\fg,m_1,m_1}^{c_m,0,s_m},\psi_K^{s_m})}\ar[rrr]^{(\psi_{\fg^*,m_1,n_1}^{1,k_1,1},\psi_K^{1})}& && \iota^*\g(n_1),K)\ar[dd]^{(1\otimes \psi_{\fg,n_1,n_1}^{c_n,0,s_n},\psi_K^{s_n})} \\
& & &\\
(\iota^*\g(m_1),K)\ar[rrr]^{(\psi_{\fg^*,m_1,n_1}^{1,k_2,1},\psi_K^{1})}&&&(\iota^*\g(n_1),K) 
}
\]
Noting that $\psi_{\fg^*,m_1,m_1}^{c_m,0,s_m}=1\otimes \psi_{\fg,m_1,m_1}^{c_m,0,s_m}$ and $\psi_{\fg^*,n_1,n_1}^{c_n,0,s_n}=1\otimes \psi_{\fg,n_1,n_1}^{c_n,0,s_n}$ 
and using Lemma \ref{lem1}, commutativity of the diagram is equivalent to 
\begin{eqnarray}\nonumber
    &&(\psi_{\fg^*,m_1,n_1}^{c_n,k_1,s_n},\psi_K^{s_n})= (\psi_{\fg^*,m_1,n_1}^{c_m,k_2,s_m},\psi_K^{s_m}).
\end{eqnarray}
Hence we must have $c_n=c_m$, $k_1=k_2$ and $s_n=s_m$. Hence the two extensions are equal.
\end{proof}
By Grothendieck's theorem any algebraic vector
bundle on $\mathbb{P}_{\C}^1$ is a direct sum of line bundles. From the last Theorem we can see what is the corresponding decomposition into line bundles of any extension over $\mathbb{P}_{\C}^1$. Here is a precise statement. 
\begin{cor}
 As algebraic $K$-equivariant vector bundles $$\left((\g(m),SO(2,\C)),  (\g(n),SO(2,\C)), (\psi_{\fg^*,m,n}^{1,k,1}, \psi^1_K)\right)$$ 
 is isomorphic to 
 \[ \mathcal{O}_{\mathbb{P}_{\C}^1}(0)\oplus \mathcal{O}_{\mathbb{P}_{\C}^1}(-k) \oplus \mathcal{O}_{\mathbb{P}_{\C}^1}(k+n-m).  \]
\end{cor}

\begin{proof}
    Since $\psi_{\fg^*,m,n}^{1,k,1}(H_m)=H_n$, the subsheaf generated by $H_m$, and  $H_n$ is isomorphic to $\mathcal{O}_{\mathbb{P}_{\C}^1}(0)$. It is the $K$ isotypic piece corresponding to the trivial representation of $K$.   Since $\psi_{\fg^*,m,n}^{1,k,1}(X_m)=x^kX_n$, the subsheaf generated by $X_m$, and  $X_n$ is isomorphic to $\mathcal{O}_{\mathbb{P}_{\C}^1}(-k)$. It is the $K$ isotypic piece corresponding to the one dimensional representation of $K$ in which $k(z)$ acts via multiplication by $z^2$. Similarly, the subsheaf generated by $Y_m$, and  $Y_n$ is isomorphic to $\mathcal{O}_{\mathbb{P}_{\C}^1}(k+n-m)$. It is the $K$ isotypic piece corresponding to the one dimensional representation of $K$ in which $k(z)$ acts via multiplication by $z^{-2}$. 
\end{proof}

\section{Declarations}
\noindent \textbf{Ethical Approval.} Not applicable.\\
\noindent \textbf{Funding.}
 This research was supported by the ISRAEL SCIENCE FOUNDATION (grant No. 1040/22).\\
\noindent \textbf{Availability of data and materials.}
Not applicable.

\bibliographystyle{plain}
\bibliography{master}
\end{document}